\documentclass[letter,10pt]{amsart}%
\usepackage{graphicx}
\usepackage{amscd}
\usepackage{amsmath}
\usepackage{amsfonts}
\usepackage{amssymb}%
\newtheorem{theorem}{Theorem}
\theoremstyle{plain}

\newtheorem{corollary}{Corollary}

\newtheorem{proposition}{Proposition}
\newtheorem{remark}{Remark}

\numberwithin{equation}{section}

\begin{document}
\author{Daniel Pellegrino}
\address{(DANIEL\ PELLEGRINO)\\
DME-Caixa Postal 10044- Campina Grande-PB-Brazil }
\email{dmp@dme.ufcg.edu.br}
\title[Scalar-valued nonlinear absolutely summing mappings]{On scalar-valued nonlinear absolutely summing mappings }
\maketitle

\begin{abstract}
In this note we investigate cases (coincidence situations) in which every
scalar-valued continuous $n$-homogeneous polynomials ($n$-linear mappings) is
absolutely $(p;q)$-summing. We extend some well known coincidence situations
and obtain several non-coincidence results, inspired in a linear technique due
to Lindenstrauss and Pe\l czy\'{n}ski$.$\ \ \ \ \ \ \ \ \ \ \ \ 

\end{abstract}

\section{Introduction}

Throughout this note $X,X_{1},...,X_{n},Y$ will stand for Banach spaces and
the scalar field $\mathbb{K}$ can be either the real or complex numbers$.$

An $m$-homogeneous polynomial $P$ from $X$ into $Y$ is said to be absolutely
$(p;q)$-summing ($p\geq\frac{q}{m}$) if there is a constant $L$ so that%

\begin{equation}
(\sum_{j=1}^{k}\Vert P(x_{j})\Vert^{p})^{\frac{1}{p}}\leq L\left\Vert
(x_{j})_{j=1}^{k}\right\Vert _{w,q}^{m}%
\end{equation}
for every natural $k$, where $\left\Vert (x_{j})_{j=1}^{k}\right\Vert
_{w,q}=\sup_{\varphi\in B_{X}%
%TCIMACRO{\U{b4}}%
%BeginExpansion
\acute{}%
%EndExpansion
}(\sum_{j=1}^{k}\mid\varphi(x_{j})\mid^{q})^{\frac{1}{q}}$. This is a natural
generalization of the concept of $(p;q)$-summing operators and in the last
years has been studied by several authors. The infimum of the $L>0$ for which
the inequality holds defines a norm $\Vert.\Vert_{as(p;q)}$ for the case
$p\geq1$ or a $p$-norm for the case $p<1$ on the space of $(p;q)$-summing
homogeneous polynomials. The space of all $m$-homogeneous $(p;q)$-summing
polynomials from $X$ into $Y$ is denoted by $\mathcal{P}_{as(p;q)}(^{m}X;Y)$
($\mathcal{P}_{as(p;q)}(^{m}X)$ if $Y=\mathbb{K}$)$.$ When $p=\frac{q}{m}$ we
have an important particular case, since in this situation there is an
analogous of the Grothendieck-Pietsch Domination Theorem. The $(\frac{q}%
{m};q)$-summing $m$-homogeneous polynomials from $X$ into $Y$ are said to be
$q$-dominated and this space is denoted by $\mathcal{P}_{d,q}(^{m}X;Y)$
($\mathcal{P}_{d,q}(^{m}X)$ if $Y=\mathbb{K}$)$.$ To denote the Banach space
of all continuous $m$-homogeneous polynomials $P$ from $X$ into $Y$ with the
$\sup$ norm we use $\mathcal{P}(^{m}X,Y)$ ($\mathcal{P}(^{m}X),$ if $Y$ is the
scalar field). Analogously, the space of all continuous $m$-linear mappings
from $X_{1}\times...\times X_{m}$ into $Y$ (with the $\sup$ norm) if
represented by $\mathcal{L}(X_{1},...,X_{m};Y)$ ($\mathcal{L}(X_{1}%
,...,X_{m})$ if $Y=\mathbb{K}$). The concept of absolutely summing multilinear
mapping follows the same pattern (for details we refer to \cite{studia}).
Henceforth every polynomial and multilinear mapping is supposed to be
continuous and every $\mathcal{L}_{p}$-space is assumed to be infinite-dimensional.

A natural problem is to find situations in which the space of absolutely
summing polynomials coincides with the space of continuous polynomials
(coincidence situations). When $Y$ is the scalar-field, these situations are
not rare as we can see on the next two well known results:

\begin{theorem}
(Matos \cite{Matos2}) \label{aa} Every scalar-valued $n$-linear mapping is
absolutely $\left(  1;1\right)  $-summing. In particular, every scalar-valued
$n$-homogeneous polynomial is absolutely $\left(  1;1\right)  $-summing (and,
a fortiori, $(q;1)$-summing for every $q\geq1$).
\end{theorem}

\begin{theorem}
\label{b}(D.P\'{e}rez-Garc\'{\i}a \cite{Perez})\label{Perezz} If $n\geq2$ and
$X$ is an $\mathcal{L}_{\infty}$-space, then every scalar-valued $n$-linear
mapping on $X$ is $(1;2)$-summing. In particular, every scalar-valued
$n$-homogeneous polynomial on $X$ is $(1;2)$-summing (and, a fortiori,
$(q;2)$-summing for every $q\geq1$).
\end{theorem}

The proof of theorem \ref{aa}, in \cite{Matos2}, is credited to A. Defant and
J. Voigt. The case $m=2$ of Theorem \ref{b} was previously proved by Botelho
\cite{Botelho} and is the unique known coincidence result for dominated
polynomials. In the Section $2$ we obtain new coincidence situations,
extending the Theorems 1 and 2. The Section $3$ has a different purpose: to
explore a technical estimate (hidden in \cite{studia}) and its several
consequences. In particular, it is shown that the Theorems 1 and 2 can not be
generalized in some other directions, and converses for the aforementioned
theorems are obtained.

\section{Coincidence situations}

The next theorem, inspired on a result of C.A. Soares, lead us to extensions
of the two theorems stated in the first section:

\begin{theorem}
\label{t1}Let $A\in\mathcal{L}(X_{1},...,X_{n};Y)$ and suppose that there
exists $C>0$ so that for any $x_{1}\in X_{1},....,x_{r}\in X_{r},$ the
$s$-linear ($s=n-r$) mapping $A_{x_{1}....x_{r}}(x_{r+1},...,x_{n}%
)=A(x_{1},...,x_{n})$ is absolutely $(1;q_{1},...,q_{s})$-summing and
$\left\Vert A_{x_{1}....x_{r}}\right\Vert _{as(1;q_{1},...,q_{s})}\leq
C\left\Vert A\right\Vert \left\Vert x_{1}\right\Vert ...\left\Vert
x_{r}\right\Vert $. Then $A$ is absolutely $(1;1,...,1,q_{1},...,q_{s})$-summing.
\end{theorem}

Proof. For $x_{1}^{(1)},...,x_{1}^{(m)}\in X_{1},....,x_{n}^{(1)}%
,...,x_{n}^{(m)}\in X_{n}$, let us consider $\varphi_{j}\in B_{Y^{\prime}}$
such that
\[
\left\Vert A(x_{1}^{(j)},...,x_{n}^{(j)})\right\Vert =\varphi_{j}%
(A(x_{1}^{(j)},...,x_{n}^{(j)}))
\]
$\text{ for every }j=1,...,m.$ Thus, defining by $r_{j}(t)$ the Rademacher
functions on $[0,1]$ and denoting by $\lambda$ the Lebesgue measure in
$I=[0,1]^{r},$ we have
\begin{align*}
&  \int\nolimits_{I}\sum\limits_{j=1}^{m}\left(  \prod_{l=1}^{r}r_{j}%
(t_{l})\right)  \varphi_{j}A(\sum\limits_{j_{1}=1}^{m}r_{j_{1}}(t_{1}%
)x_{1}^{(j_{1})},...,\sum\limits_{j_{r}=1}^{m}r_{j_{r}}(t_{r})x_{r}^{(j_{r}%
)},x_{r+1}^{(j)},...,x_{n}^{(j)})d\lambda\\
&  =\sum\limits_{j,j_{1},...j_{r}=1}^{m}\varphi_{j}A(x_{1}^{(j_{1})}%
,...,x_{r}^{(j_{r})},x_{r+1}^{(j)},...,x_{n}^{(j)})\int\limits_{0}^{1}%
r_{j}(t_{1})r_{j_{1}}(t_{1})dt_{1}...\int\limits_{0}^{1}r_{j}(t_{r})r_{j_{r}%
}(t_{r})dt_{r}\\
&  =\sum\limits_{j=1}^{m}\sum\limits_{j_{1}=1}^{m}...\sum\limits_{j_{r}=1}%
^{m}\varphi_{j}A(x_{1}^{(j_{1})},...,x_{r}^{(j_{r})},x_{r+1}^{(j)}%
,...,x_{n}^{(j)})\delta_{jj_{1}}...\delta_{jj_{r}}\\
&  =\sum\limits_{j=1}^{m}\varphi_{j}A(x_{1}^{(j)},...,x_{n}^{(j)}%
)=\sum\limits_{j=1}^{m}\left\Vert A(x_{1}^{(j)},...,x_{n}^{(j)})\right\Vert
=(\ast).
\end{align*}
So, for each $l=1,...,r,$ assuming $z_{l}=\sum\limits_{j=1}^{m}r_{j}%
(t_{l})x_{l}^{(j)}$ we obtain
\begin{align*}
(\ast)  &  =\int\nolimits_{I}\sum\limits_{j=1}^{m}\left(  \prod_{l=1}^{r}%
r_{j}(t_{l})\right)  \varphi_{j}A(\sum\limits_{j_{1}=1}^{m}r_{j_{1}}%
(t_{1})x_{1}^{(j_{1})},...,\sum\limits_{j_{r}=1}^{m}r_{j_{r}}(t_{r}%
)x_{r}^{(j_{r})},x_{r+1}^{(j)},...,x_{n}^{(j)})d\lambda\\
&  \leq\int\nolimits_{I}\left\vert \sum\limits_{j=1}^{m}\left(  \prod
_{l=1}^{r}r_{j}(t_{l})\right)  \varphi_{j}A(\sum\limits_{j_{1}=1}^{m}r_{j_{1}%
}(t_{1})x_{1}^{(j_{1})},...,\sum\limits_{j_{r}=1}^{m}r_{j_{r}}(t_{r}%
)x_{r}^{(j_{r})},x_{r+1}^{(j)},...,x_{n}^{(j)})\right\vert d\lambda\\
&  \leq\int\nolimits_{I}\sum\limits_{j=1}^{m}\left\Vert A(\sum\limits_{j_{1}%
=1}^{m}r_{j_{1}}(t_{1})x_{1}^{(j_{1})},...,\sum\limits_{j_{r}=1}^{m}r_{j_{r}%
}(t_{r})x_{r}^{(j_{r})},x_{r+1}^{(j)},...,x_{n}^{(j)})\right\Vert d\lambda\\
&  \leq\sup_{t_{l}\in\lbrack0,1],l=1,...,r}\sum\limits_{j=1}^{m}\left\Vert
A(\sum\limits_{j_{1}=1}^{m}r_{j_{1}}(t_{1})x_{1}^{(j_{1})},...,\sum
\limits_{j_{r}=1}^{m}r_{j_{r}}(t_{r})x_{r}^{(j_{r})},x_{r+1}^{(j)}%
,...,x_{n}^{(j)})\right\Vert \\
&  \leq\sup_{t_{l}\in\lbrack0,1],l=1,...,r}\left\Vert A_{z_{1}...z_{r}%
}\right\Vert _{as(1;q_{1},...,q_{s})}\left\Vert (x_{r+1}^{(j)})_{j=1}%
^{m}\right\Vert _{w,q_{1}}...\left\Vert (x_{n}^{(j)})_{j=1}^{m}\right\Vert
_{w,q_{s}}\\
&  \leq\sup_{t_{l}\in\lbrack0,1],l=1,...,r}C\left\Vert A\right\Vert \left\Vert
z_{1}\right\Vert ...\left\Vert z_{r}\right\Vert \left\Vert (x_{r+1}%
^{(j)})_{j=1}^{m}\right\Vert _{w,q_{1}}...\left\Vert (x_{n}^{(j)})_{j=1}%
^{m}\right\Vert _{w,q_{s}}\\
&  \leq C\left\Vert A\right\Vert \left(  \prod_{l=1}^{r}\left\Vert
(x_{l}^{(j)})_{j=1}^{m}\right\Vert _{w,1}\right)  \left(  \prod_{l=r+1}%
^{n}\left\Vert (x_{l}^{(j)})_{j=1}^{m}\right\Vert _{w,q_{l}}\right)  .\text{ }%
\end{align*}
We have the following straightforward consequence, generalizing Theorem
\ref{aa}:

\begin{corollary}
\label{ssss}If
\[
\mathcal{L}(X_{1},...,X_{m};Y)=\mathcal{L}_{as(1;q_{1},...,q_{m})}%
(X_{1},...,X_{m};Y)
\]
then, for any Banach spaces $X_{m+1},...,X_{n},$ we have
\[
\mathcal{L}(X_{1},...,X_{n};Y)=\mathcal{L}_{as(1;q_{1},...,q_{m}%
,1,...,1)}(X_{1},...,X_{n};Y).
\]

\end{corollary}

The following corollary (whose proof is simple and we omit) is consequence of
the Theorems \ref{b} and \ref{t1}.

\begin{corollary}
\label{c1}If $X_{1},...,X_{s}$ are $\mathcal{L}_{\infty}$-spaces then, for any
choice of Banach spaces $X_{s+1},...,X_{n},$ we have
\[
\mathcal{L}(X_{1},...,X_{n})=\mathcal{L}_{as(1;q_{1},...,q_{n})}%
(X_{1},...,X_{n}),
\]
where $q_{1}=...=q_{s}=2$ and $q_{s+1}=....=q_{n}=1.$
\end{corollary}

It is obvious that Corollary \ref{c1} is still true if we replace the scalar
field by any finite dimensional Banach space. A natural question is whether
Corollary \ref{c1} can be improved for some infinite dimensional Banach space
in the place of $\mathbb{K}.$ Precisely, the question is:

\begin{itemize}
\item If $X_{1},...,X_{k}$ are $\mathcal{L}_{\infty}$-spaces, is there some
infinite dimensional Banach space $Y$ such that
\[
\mathcal{L}(X_{1},...,X_{k},...,X_{n};Y)=\mathcal{L}_{as(1;q_{1},....,q_{n}%
)}(X_{1},...,X_{k},...,X_{n};Y),
\]
where $q_{1}=...=q_{k}=2$ and $q_{k+1}=....=q_{n}=1,$ regardless of the Banach
spaces $X_{k+1},...,X_{n}$?
\end{itemize}

The answer to this question is no, as we can see on the following proposition:

\begin{proposition}
Suppose that $X_{1},...,X_{k}$ are infinite dimensional $\mathcal{L}_{\infty}%
$-spaces. If $q_{1}=...=q_{k}=2$, $q_{k+1}=....=q_{n}=1$ and
\[
\mathcal{L}(X_{1},...,X_{k},...,X_{n};Y)=\mathcal{L}_{as(1;q_{1},....,q_{n}%
)}(X_{1},...,X_{k},...,X_{n};Y),
\]
regardless of the Banach spaces $X_{k+1},...,X_{n},$ then $\dim Y<\infty.$
\end{proposition}

Proof. By a standard localization argument, it suffices to prove that if $\dim
Y=\infty,$ then
\[
\mathcal{L}(^{n}c_{0};Y)\neq\mathcal{L}_{as(1;q_{1},....,q_{n})}(^{n}%
c_{0};Y),
\]
where $q_{1}=...=q_{k}=2$ and $q_{k+1}=....=q_{n}=1.$ But, from{
\ }\cite[Theorem 8]{studia} we have%
\[
\mathcal{L}(^{n}c_{0};Y)\neq\mathcal{L}_{as(q;q_{1},....,q_{n})}(^{n}%
c_{0};Y),
\]
regardless of the $q<2$ and $q_{1},...,q_{n}\geq1.$

\section{ Non-coincidence situations}

\bigskip Assume that $X$ is an infinite dimensional Banach space and suppose
that $X$ has a normalized unconditional Schauder basis $(x_{n})$ with
coefficient functionals $(x_{n}^{\ast}).$ If $\mathcal{P}_{as(q;1)}%
(^{m}X;Y)=\mathcal{P}(^{m}X;Y),$ it is natural to ask:

What is the best $t$ such that in this situation $(x_{n}^{\ast}(x))\in l_{t}$
for each $x\in X$ ? The best such $t$ will be denoted by $\mu=\mu(X,Y,q,m).$

In \cite{studia}, inspired by a linear result due to Lindenstrauss and
Pe\l czy\'{n}ski, we have proved:

\begin{theorem}
(Pellegrino \cite[Theorem 5]{studia}) Let $X$ and $Y$ be infinite dimensional
Banach spaces. Suppose that $X$ has an unconditional Schauder basis $(x_{n}%
)$.\ If $Y$ finitely factors the formal inclusion $l_{p}\rightarrow l_{\infty
}$ and $\mathcal{P}_{as(q;1)}(^{m}X;Y)=\mathcal{P}(^{m}X;Y)$ with $\frac{1}%
{m}\leq q,$ then

(a) $\mu\leq\frac{mpq}{p-q}$ if $q<p$

(b) $\mu\leq mq$ if $q\leq\frac{p}{2}.$
\end{theorem}

However, observing the proof of this theorem in \cite[Theorem 5]{studia}, one
can see that it is absolutely not necessary to assume that $\dim Y=\infty.$
Only in Corollary 1 of \cite{studia} (when Dvoretzky-Rogers Theorem is
invoked) it is indeed necessary to assume $dimY$$=$$\infty$. A slight change
on the proof of \cite[Theorem 5]{studia} yields:

\begin{theorem}
\label{main}\label{azz} Let $X$ be an infinite dimensional Banach space with a
normalized unconditional Schauder basis $(x_{n})$.\ If $\mathcal{P}%
_{as(q;1)}(^{m}X)=\mathcal{P}(^{m}X),$ then

(a) $\mu\leq\frac{mq}{1-q}$ if $q<1.$

(b) $\mu\leq mq$ if $q\leq\frac{1}{2}.$
\end{theorem}

Proof. If $x=\sum\limits_{j=1}^{\infty}a_{j}x_{j}$ and $\{\mu_{i}\}_{i=1}^{n}$
is such that $\sum\limits_{j=1}^{n}\mid\mu_{j}\mid^{\frac{1}{q}}=1,$ define
$P:X\rightarrow\mathbb{K}$ by $Px=\sum\limits_{j=1}^{n}\left\vert \mu
_{j}\right\vert ^{\frac{1}{q}}a_{j}^{m}$ and proceed as in the proof \ of
Theorem 5 of \cite{studia}, with $p=1$.

\bigskip Now we list several important consequences of Theorem \ref{azz}. For
example, the Corollaries \ref{mm} and \ref{nn} give converses for Theorems 1
and 2, respectively. The proof of Corollaries \ref{mm}, \ref{nn}, \ref{s} and
\ref{d} are immediate (using Theorem \ref{azz} and standard localizations
techniques in order to extend the results from $c_{0}$ to $\mathcal{L}%
_{\infty}$-spaces):

\begin{corollary}
\bigskip\label{mm}Let $m$ be a fixed natural number. Then $\mathcal{P}%
_{as(q;1)}(^{m}X)=\mathcal{P}(^{m}X)$ for every $X$ if and only if $q\geq1.$
\end{corollary}

\begin{corollary}
\label{nn}If $m\geq2$ and $X$ is an $\mathcal{L}_{\infty}$-space, then
$\mathcal{P}_{as(q;2)}(^{m}X)=\mathcal{P}(^{m}X)$ if and only if $q\geq1.$
\end{corollary}

\begin{corollary}
\label{s}If $X$ is an $\mathcal{L}_{\infty}$-space, then $\mathcal{P}%
_{d,q}(^{m}X)\neq\mathcal{P}(^{m}X)$ for every $q<m.$
\end{corollary}

\begin{corollary}
\label{d}If $X$ is an $\mathcal{L}_{\infty}$-space, then $\mathcal{P}%
_{d,q}(^{2}X)=\mathcal{P}(^{2}X)$ if and only if $q\geq2.$
\end{corollary}

We also have:

\begin{corollary}
If $q\leq\frac{1}{2}$ and $X$ is an $\mathcal{L}_{p}$ space ($p\geq1$), then
$\mathcal{P}_{as(q;1)}(^{m}X)=\mathcal{P}(^{m}X)$ if and only if $p\leq mq.$
\end{corollary}

\bigskip Proof. If \bigskip$\mathcal{P}_{as(q;1)}(^{m}X)=\mathcal{P}(^{m}X),$
the Theorem \ref{main} assures that $p\leq mq.$ On the other hand, if $p\leq
mq$ and $P\in\mathcal{P}(^{m}X),$ then
\begin{align*}
(\sum_{j=1}^{k}\Vert P(x_{j})\Vert^{q})^{\frac{1}{q}}  &  \leq\Vert
P\Vert(\sum_{j=1}^{k}\Vert x_{j}\Vert^{mq})^{\frac{1}{q}}\\
&  \leq\Vert P\Vert(\sum_{j=1}^{k}\Vert x_{j}\Vert^{p})^{\frac{m}{p}}\\
&  \leq\Vert P\Vert\left\Vert (x_{j})_{j=1}^{k}\right\Vert _{w,1}^{m},
\end{align*}
where the last inequality holds since $l_{p}$ has cotype $p$ (for $p\geq1$)
and thus $id:l_{p}\rightarrow l_{p}$ is absolutely $(p;1)$-summing.

All these results can be adapted (including Theorem \ref{main}), mutatis
mutandis, to the multilinear cases. In particular, one can extend Corollary
\ref{c1}:

\begin{corollary}
\label{swq}Let $X_{1},...,X_{s}$ be $\mathcal{L}_{\infty}$-spaces,
$q_{1}=...=q_{s}=2$ and $q_{s+1}=....=q_{n}=1$. We have
\[
\mathcal{L}(X_{1},...,X_{n})=\mathcal{L}_{as(q;q_{1},...,q_{n})}%
(X_{1},...,X_{n}),
\]
for any choice of Banach spaces $X_{s+1},...,X_{n},$ if and only if $q\geq1.$
\end{corollary}

\begin{remark}
For the bilinear case it is not hard to prove that when $X$ is an
$\mathcal{L}_{\infty}$-space, $\mathcal{L}_{d,q}(^{2}X)=\mathcal{L}(^{2}X)$ if
and only if $q\geq2.$ However, this result can not be straightforwardly
adapted for polynomials. Non-coincidence results for multilinear mappings in
general does not imply non-coincidence results for polynomials.
\end{remark}

\end{document}